\documentclass[11pt]{amsart}

\usepackage[utf8]{inputenc}
\usepackage[T1]{fontenc}
\usepackage[french, english]{babel}

\usepackage{amsthm}
\usepackage{amsmath}
\usepackage{amssymb}
\usepackage{mathrsfs}
\usepackage{bbm}
\usepackage{graphicx}
\usepackage{listings}
\usepackage{mathrsfs}
\usepackage{enumerate}
\usepackage{pict2e}
\usepackage{stmaryrd}
\usepackage{mathtools}

\usepackage[all,cmtip]{xy}

\usepackage{multirow}

\usepackage{color}
\definecolor{marin}{rgb}   {0.,   0.1,   0.5} 
\definecolor{rouge}{rgb}   {0.8,   0.,   0.} 
\definecolor{sepia}{rgb}   {0.4,   0.25,   0.} 
\definecolor{mag}{rgb}   {0.3,   0,   0.3} 

\usepackage[colorlinks,citecolor=sepia,linkcolor=marin,urlcolor=mag ,bookmarksopen,bookmarksnumbered]{hyperref}

\newtheorem{theorem}{Theorem}[section]
\newtheorem{corollary}[theorem]{Corollary}

\newtheorem{definition}[theorem]{Definition}

\newcommand{\ic}{\mathrm{i}}

\newcommand{\dd}{\mathrm{d}}

\newcommand{\Z}{\mathbb{Z}}

\newcommand{\N}{\mathbb{N}}

\newcommand{\R}{\mathbb{R}}

\newcommand{\eps}{\varepsilon}

       \renewcommand{\S}{{\mathcal S}}

\begin{document}
\title[On almost periodic solutions to NLS without external parameters]{On almost periodic solutions to NLS without external parameters}

\author{Joackim Bernier}

\address{\small{Nantes Universit\'e, CNRS, Laboratoire de Math\'ematiques Jean Leray, LMJL,
F-44000 Nantes, France
}}

\email{joackim.bernier@univ-nantes.fr}

\author{Beno\^it Gr\'ebert}

\address{\small{Nantes Universit\'e, CNRS, Laboratoire de Math\'ematiques Jean Leray, LMJL,
F-44000 Nantes, France
}}

\email{benoit.grebert@univ-nantes.fr}



\begin{abstract} In this note, we present a result established  in \cite{BGR24} where we prove that nonlinear Schr\"odinger equations on the circle, without external parameters, admit plenty of infinite dimensional non resonant invariant tori, or equivalently, plenty of almost periodic solutions. Our aim is to propose  an extended sketch of the proof, emphasizing the new points which have enabled us to achieve this result.
\end{abstract} 
\maketitle


\section{Result and context}

 We consider the  nonlinear Schr\"odinger equation on the circle $\mathbb T:=\R/2\pi\Z$
\begin{equation}
\tag{NLS}
\label{eq:NLS}
\ic \partial_t u + \partial_x^2 u = f(|u|^2)u, \quad x\in \mathbb{T},\ t\in\R
\end{equation}
where $f$ is a real entire function such that $f(0)=0$, $f'(0)\neq0$  (i.e. the cubic term is present) and $z\mapsto f(z^2)$ grows at most exponentially fast\footnote{i.e. there exists $C>0$ such that for all $z\in \mathbb{C}$, $|f(z)|\leq C \exp(C\sqrt{|z|})$.}. We point out that \eqref{eq:NLS} is a classical and emblematic Hamiltonian system. 
When considering a Hamiltonian system, we are naturally led to look for invariant sets.    An essential challenge  is to construct invariant sets that are "large" in the sense that they occupy a significant part of the infinite-dimensional phase space.

\medskip

We are looking for solutions with Sobolev regularity in the space variable\footnote{in \cite{BGR24} we also consider analytic solutions but, for simplicity, we prefer to choose $H^s$ as phase space in this note.}, i.e. we consider solutions living in the phase space
$$
H^s:= \big\{ u \in L^2(\mathbb{T}) \ | \ \| u\|_{H^s}^2 := \sum_{k\in \mathbb{Z}} |u_k|^2 \langle k\rangle^{2s}  < \infty\big \}
$$
where we always identify functions $ u \in L^2(\mathbb{T})$ with their Fourier coefficients 
\begin{equation}
\label{eq:def_Fourier}
 u_k:=(2\pi)^{-1}\int_{\mathbb{T}}u(x)e^{ikx}\dd x, \quad k\in \mathbb{Z}.
 \end{equation}
The phase space $H^s$ is foliated by invariant tori for the linear part of the equation \eqref{eq:NLS}. Namely
for all $\xi\in (\mathbb{R}_+)^\mathbb{Z}$ the torus
\begin{equation}
\label{eq:def_T_droit}
\mathrm T_\xi :=\{ u\in H^s\mid \forall k\in\mathbb{Z}, \quad  |u_k|^2=\xi_k\, \},
\end{equation}
is an invariant set for the system
\begin{equation}
\label{eq:LS}
\ic \partial_t u +\partial_x^2 u=0 .
\end{equation}
The natural question is:
\begin{center}
\emph {Which linear tori persist, slightly deformed, as invariant sets of \eqref{eq:NLS}?}
 \end{center}
 Kuksin--P\"oschel  proved in \cite{KP96} that most of the finite dimensional tori, close to the origin, persist as KAM tori (i.e. Kronecker tori that are linearly stable). Unfortunately, their result suffers a limitation : the larger the dimension of the torus is, the closer to the origin it has to be. This limitation is not specific to the work of Kuksin--P\"oschel, it appears in all the papers about existence of finite dimensional invariant tori. Thus, up to now, there have been no results concerning the existence of infinite dimensional invariant tori for  \eqref{eq:NLS} and, more generally, for (non-integrable) Hamiltonian PDEs without external parameter. Let us define the kind of tori we are going to prove the persistency of:

\begin{definition}[Non resonant infinite dimensional Kronecker tori] A subset $\mathcal{T}$ of $H^s$ is a \emph{non resonant infinite dimensional Kronecker torus} (for \eqref{eq:NLS}) if there exists a homeomorphism\footnote{Naturally, we always equip $\mathbb{T}^\mathbb{N}$ with the product topology.} $\Psi : \mathbb{T}^\mathbb{N} \to \mathcal{T}$ and a sequence of rationally independent numbers $\omega \in \mathbb{R}^{\mathbb{N}}$ such that
 \begin{center}
 $\forall \theta\in \mathbb{T}^{\mathbb{N}}$, $t\mapsto \Psi( \omega t +\theta)$ is a global  solution to  \eqref{eq:NLS}.
\end{center}
\end{definition}
Since the frequencies are rationally independent (or non resonant), i.e. $k\cdot\omega\neq 0$ for all $0\neq k\in\Z^\N$ of finite support, on these tori each orbit is dense and the flow is ergodic (for the product measure). Our main result can be summarized in a very short theorem:
\begin{theorem} \label{thm:0} There exists non resonant infinite dimensional Kronecker tori for \eqref{eq:NLS}.
\end{theorem}
Actually, as explain in section \ref{sec:scheme}, we construct a whole family of non resonant infinite dimensional Kronecker tori which accumulates on the finite dimensional  tori constructed by Kuksin--P\"oschel in \cite{KP96}.

\medskip

Another classical way of considering this kind of problem, consists in looking for almost periodic solution (which are not quasi periodic). Up to now, there have been no results concerning the existence of such solutions for  \eqref{eq:NLS} and more generally for not completely integrable Hamiltonian PDEs without external parameter.
\begin{definition}[Quasi-periodic  function] A continuous function $q:\R\to\R$ is called quasi-periodic with frequencies $\omega=(\omega_1,\cdots,\omega_n)$, if there exists a continuous function, $Q:\mathbb T^n\to \R$, such that $q(t)=Q(\omega t)$ for all $t\in\R$.
\end{definition}
Then following H. Bohr (see \cite{Boh47}) we  recall the definition of an almost-periodic function:
\begin{definition}[Almost-periodic function] A continuous function  is called almost-periodic if it is a uniform limit of quasi-periodic functions.
\end{definition}
As a corollary of Theorem \ref{thm:0}, we get:
\begin{corollary}
There exist almost periodic solutions to \eqref{eq:NLS} which are not quasi periodic.
\end{corollary}

 \noindent \underline{ About the related literature:} For several decades now, many mathematicians have been interested in KAM theory for nonlinear Hamiltonian PDEs. For finite dimensional Hamiltonian system, standard KAM theory ensures that, under some non-degeneracy assumptions, most tori survive and thus most of the solutions are quasi-periodic (see e.g. \cite{KP03}). 
 Naturally, in an infinite-dimensional phase space, infinite-dimensional invariant sets are expected. In other words, for PDEs, \emph{we expect almost periodic solutions to be typical}.

Many results have been obtained, but mainly concerning the existence of finite-dimensional invariant tori and often at the cost of adding external parameters to the equation to simplify the treatment of the small denominators inherent in KAM theory (see e.g. \cite{Kuk87, Way90, CW93, Bou94, Pos96, KP96} for the earliest results for 1d PDEs, \cite{Bou98, EK10, BB13, PP15,  EGK16, GP16} for some results in higher dimension, \cite{BBHM18, FG24, BHM23} for some results about quasi-linear equations and \cite{BKM18} for large amplitude solutions). Most of these results allow to get invariant tori of arbitrarily large dimension. Nevertheless, in the proofs, the larger the dimension is, the smaller the perturbation parameter has to be. Therefore, they say nothing about the existence of infinite dimensional invariant tori.

\medskip

  The first result in a PDE context is due to Bourgain in \cite{Bou96} who proved the existence of almost periodic solutions to the wave equation in presence of a random potential. Then, few years later, in \cite{Bou05}, the same author considered \eqref{eq:NLS} but with an extra random convolution potential $V$ as external parameter:
\begin{equation}
\label{eq:NLSconv}
\ic \partial_t u + \partial_x^2 u+ V*u = f(|u|^2)u, \quad x\in \mathbb{T},\ t\in\R.
\end{equation}
Again he proved the existence of almost periodic solutions but, this time, with more reasonable decrease  of the Fourier coefficients of the solution. This result has been revisited and improved by a series of paper (in particular \cite{BMP20,Con23} and  
\cite{CY21} for the wave equation). Further, using sparsity arguments, Biasco--Masetti--Procesi proved in \cite{BMP24} the existence of almost periodic solutions of finite regularity. In all these results the Fourier coefficients of the potential $V$ are bounded and does not converge to $0$. This is a critical ingredient of these proofs which is used to deal with the small divisors (since the linear frequencies for \eqref{eq:NLSconv}  are $\omega_j=j^2-V_j,\ j\in\Z$). 
This assumption on $V$ seems to be the major obstruction  to remove $V$ in \eqref{eq:NLSconv} (and so to get \eqref{eq:NLS}).  Indeed, following an observation due to  Kuksin--P\"oschel in \cite{KP96},  \eqref{eq:NLS} is somehow equivalent to \eqref{eq:NLSconv} where $V$ would depend nonlinearly on $u$ through the relation $V_j = |u_j|^2$, $j\in \mathbb{Z}$. It therefore makes sense to try to match the regularity of $V$ and $u$ in \eqref{eq:NLSconv} in order to get rid of the external parameters. In this direction, introducing a new approach based on tree expansions and renormalization groups, Corsi--Gentile--Procesi (see \cite{CGP24})  succeeded in greatly improving the regularity of $V$, indeed they can consider $V$ of any finite regularity $C^n$. Unfortunately, their almost periodic solutions are still Gevrey in the space variable. Very recently, Biasco--Corsi--Gentile--Procesi extended this result (in \cite{BCGP24}) by proving that their set of almost periodic solutions is asymptotically of full measure.

\medskip

Without external parameters, the only known results prove that most of the linear infinite dimensional tori are almost preserved for very long times (see e.g. \cite{Bou00,BFG20,LX24,BC24} for \eqref{eq:NLS} and \cite{BG21} for perturbations of KdV or Benjamin--Ono equations). The proofs of these results, somewhere between Birkhoff normal forms and KAM, deal with many similar difficulties (in particular the degeneracy of small divisors) but make crucial use of the fact that time, although potentially very long, is finite.

  \medskip

 \noindent \underline{ About our approach:} Our approach is different from that adopted in the above-mentioned articles which concern infinite-dimensional tori (in particular in \cite{Bou05, BMP20}): we don't try to construct an infinite-dimensional torus directly, but, iteratively,
a convergent sequence of invariant tori whose dimension goes to $+\infty$.
 This was the scheme of proof proposed by P\"oschel in \cite{Pos02} but then rarely used because it requires the nonlinearity to be regularizing which is not the case for \eqref{eq:NLS} neither for most of the Hamiltonian PDEs (see also \cite{GX13}). Clearly this is why we need the regularization normal form result described in section \ref{sec:reg}.
 However it is not the end of the story since, first, the regularization is just partial, but also P\"oschel was considering external parameters (namely the Fourier coefficients of a potential) and we don't. 
 
 \medskip

 \noindent  In the next two sections we will give the principal keys of the proof. The two sections are independent, the result of section \ref{sec:reg}, the regularizing normal form, will be used as a black box in section \ref{sec:scheme} where we will focus on the principal difficulty: how to deal with internal parameters when constructing an infinite dimensional invariant torus.

\section{About the regularizing normal form}\label{sec:reg} The objective of this section is to provide the intuitions that led us to obtain our regularizing normal form result. First, we discuss the regularizing effects given by the Birkhoff normal form. These effects have numerous limitations. In particular, we observe that, at best, we can hope that the nonlinearity becomes smoothing modulo a gauge transform. Secondly, we explain how to adapt the Birkhoff normal form procedure to achieve such a regularizing effect. For more results about regularizing effects generated by the dispersion we refer to \cite{ErTz1,ErTz3,McConnell}.

\medskip

Before starting the explanations, let us briefly recall the structure of \eqref{eq:NLS}. Indeed, \eqref{eq:NLS} is a Hamiltonian system: it rewrites
$$
\ic \partial_t u = \nabla H(u)
$$
where $\nabla$ denotes the $L^2$ gradient,
$$
H(u) = \underbrace{\frac12 \sum_{k\in \mathbb{Z}} k^2 u_k^2 }_{=:Z_2(u)}+ \sum_{q\geq 2} a_{2q} \| u\|_{L^{2q}}^{2q} \quad \mathrm{with} \quad a_{2q} = \frac{f^{(q-1)}(0)}{2(q!)}.
$$
It is important to note that
$$
\| u\|_{L^{2q}}^{2q} =\frac{(q!)^2}{(2q)!} \sum_{\substack{\boldsymbol{\sigma}_1 + \cdots +\boldsymbol{\sigma}_{2q} = 0 \\ \boldsymbol{\sigma}_1 \boldsymbol{\ell}_1+ \cdots +\boldsymbol{\sigma}_{2q}  \boldsymbol{\ell}_{2q}} =0 } u_{\boldsymbol{\ell}_1}^{\boldsymbol{\sigma}_1} \cdots  u_{\boldsymbol{\ell}_{2q}}^{\boldsymbol{\sigma}_{2q}} 
$$
where $u_k$, $k\in \mathbb{Z}$, denotes the $k$-th Fourier coefficients of $u$ and by convention $u_k^{-1} = \overline{u_k}$.

We also introduce a notation: let $\boldsymbol{\ell}\in \Z^q$,    $\boldsymbol{\ell}^*$ denotes a rearrangement of $\boldsymbol{\ell}$ such that $$|\boldsymbol{\ell}^*_1|\geq \cdots \geq |\boldsymbol{\ell}^*_{q}|.$$

\subsection{Birkhoff normal form and natural obstructions} \label{sub:Birk} As a starting point, let us recall some smoothing effects we get by putting the equation in Birkhoff normal form.

\subsubsection*{Normal form and decomposition} We fix $r\geq 2$ and we put \eqref{eq:NLS} in Birkhoff normal form up to order $2r$. Applying standard technics, we get a symplectic change of variables $\tau$, close to the identity, defined on a neighborhood of the origin in $H^s$ such that
$$
H \circ \tau^{-1} \mathop{=}_{u\to 0} Z_2 +\underbrace{ \sum_{q=2}^{r} \sum_{\substack{\boldsymbol{\sigma}_1 + \cdots +\boldsymbol{\sigma}_{2q} = 0 \\ \boldsymbol{\sigma}_1 \boldsymbol{\ell}_1+ \cdots +\boldsymbol{\sigma}_{2q}  \boldsymbol{\ell}_{2q} =0 \\ \boldsymbol{\sigma}_1 \boldsymbol{\ell}_1^2+ \cdots +\boldsymbol{\sigma}_{2q}  \boldsymbol{\ell}_{2q}^2 = 0}  } B^{\boldsymbol{\ell},\boldsymbol{\sigma}} u_{\boldsymbol{\ell}_1}^{\boldsymbol{\sigma}_1} \cdots  u_{\boldsymbol{\ell}_{2q}}^{\boldsymbol{\sigma}_{2q}}  }_{=: B}+ \mathcal{O}(u^{2r+2})
$$
where the coefficients $ B^{\boldsymbol{\ell},\boldsymbol{\sigma}} $ are real and uniformly bounded.

\medskip 

Then we split $B$ in two parts
$$
B = G + M
$$
where $M$ contain the terms of $B$ of the form $|u_a|^2 u_{\boldsymbol{\ell}_1}^{\boldsymbol{\sigma}_1} \cdots  u_{\boldsymbol{\ell}_{2p-2}}^{\boldsymbol{\sigma}_{2q-2}} $ with $a\geq |\boldsymbol{\ell}_1^*|$. More precisely, the coefficients of $M$ satisfy
$$
M^{\boldsymbol{\ell},\boldsymbol{\sigma}} = \mathbbm{1}_{(\boldsymbol{\ell},\boldsymbol{\sigma})\in \mathbf{P}_q} B^{\boldsymbol{\ell},\boldsymbol{\sigma}}
$$
where $\mathbf{P}_{q}$ denotes the set of indices $(\boldsymbol{\ell},\boldsymbol{\sigma}) \in\mathbb{Z}^{q} \times \{-1,1\}^q$ such that there is a pairing between the two highest indices: there exist $j_1,j_2\in\{1,\dots,q\}$ such that 
$$
j_1 \neq j_2, \quad  \boldsymbol{\ell}_{j_1}=\boldsymbol{\ell}_{j_2}, \quad \boldsymbol{\sigma}_{j_1}=-\boldsymbol{\sigma}_{j_2}, \quad |\boldsymbol{\ell}_{j_1}| = |\boldsymbol{\ell}_{1}^*|  .
$$

\subsubsection*{Smoothing terms} We note that $G$ is $(s-s_0)/2$ smoothing where $s_0\in(1/2,s)$ can be chosen arbitrarily close to $1/2$. More precisely, $\nabla G$ is a smooth polynomial from $H^s$ to $H^{s+ (s-s_0)/2 }$. To prove it, it suffices to note that 
\begin{equation}
\label{eq:ca_implique}
\left. \begin{array}{lll} 
\substack{\boldsymbol{\sigma}_1 + \cdots +\boldsymbol{\sigma}_{2q} = 0 \\ 
\boldsymbol{\sigma}_1 \boldsymbol{\ell}_1+ \cdots +\boldsymbol{\sigma}_{2q}  \boldsymbol{\ell}_{2q} =0 \\
 \boldsymbol{\sigma}_1 \boldsymbol{\ell}_1^2+ \cdots +\boldsymbol{\sigma}_{2q}  \boldsymbol{\ell}_{2q}^2 = 0} \\ 
(\boldsymbol{\ell},\boldsymbol{\sigma}) \notin \mathbf{P}_q  \end{array} \right\} \implies |\boldsymbol{\ell}_1^*| \lesssim_q |\boldsymbol{\ell}_3^*|^2.
\end{equation}\\
We note that, from the point of view of the coefficients, it implies that 
$$
\forall \delta \geq 0, \quad |G^{\boldsymbol{\ell},\boldsymbol{\sigma}}| \lesssim_{q,\delta} 1 \lesssim_{q,\delta} \frac{ |\boldsymbol{\ell}_3^*|^{2\delta}}{ |\boldsymbol{\ell}_1^*|^{\delta}}.
$$
{This leads to a $\delta$ smoothing effect  if you can control  $|\boldsymbol{\ell}_3^*|^{2\delta}$, i.e. for $\delta\leq (s-s_0)/2.$ }

\subsubsection*{Non-smoothing terms} There are two kind of non-smoothing terms. On the one hand, although the terms in $M$ are very special, some of them are not smoothing. For example, $B$ contains non-smoothing terms of the form
\begin{equation}
\label{eq:relou}
\| u\|_{L^2}^4 \quad  \mathrm{or} \quad \| u\|_{L^2}^2 \sum_{k\in \mathbb{Z}} |u_k|^4.
\end{equation}
These terms directly come from $H$: since they are resonant, they are not affected by the Birkhoff normal form procedure. 

\medskip

On the other hand, the Birkhoff normal form procedure is known to be divergent (the coefficients go to $+\infty$ as $r$ goes to $+\infty$) and, $r$ being fixed, the terms in $\mathcal{O}(u^{2r+2})$ have no particular structure and are not smoothing.

\subsection{Reasonable smoothing effects} \label{sub:desc_reg} The non smoothing terms \eqref{eq:relou} are integrable and so it is hopeless to remove them by any normal form procedure. Nevertheless, we note that they are very special. The fact that they are not smoothing only comes from the $\| u\|_{L^2}^2$ factors. Since the $L^2$ norm is a constant of the motion of \eqref{eq:NLS}, it is tempting to consider the $\| u\|_{L^2}^2$ factors as constants. This is exactly what we are doing but, of course, there is a price to pay.  

\medskip

Indeed, we are going to conjugate $H$ to a Hamiltonian of the form
$$
H\circ \tau^{-1} = Z_2 + P = Z_2 +  \sum_{q=2}^{\infty} \sum_{\substack{\boldsymbol{\sigma}_1 + \cdots +\boldsymbol{\sigma}_{2q} = 0 \\ \boldsymbol{\sigma}_1 \boldsymbol{\ell}_1+ \cdots +\boldsymbol{\sigma}_{2q}  \boldsymbol{\ell}_{2q} =0 }  } P^{\boldsymbol{\ell},\boldsymbol{\sigma}}(\|u\|_{L^2}^2) u_{\boldsymbol{\ell}_1}^{\boldsymbol{\sigma}_1} \cdots  u_{\boldsymbol{\ell}_{2q}}^{\boldsymbol{\sigma}_{2q}} 
$$
where, locally uniformly in $L^2$,
\begin{equation}
\label{eq:coucou}
 |P^{\boldsymbol{\ell},\boldsymbol{\sigma}}(\|u\|_{L^2}^2)| \lesssim_{\delta,q} 1\wedge \frac{ |\boldsymbol{\ell}_3^*|^{2\delta}}{ |\boldsymbol{\ell}_1^*|^{\delta}}.
\end{equation}
with $\delta=1$. Such a normal form does not mean that we have conjugated \eqref{eq:NLS} to a Hamiltonian system with a smoothing nonlinearity. Indeed, the vector field of $P$ is of the form
$$
\nabla P = 2 (\partial_\mu P(u) ) u + 2 \nabla_{v} P(u)
$$
where $\mu$ refers to the dependency of the coefficients with respect to the $L^2$ norm of $u$ (and $v$ to the other factors). Thanks to \eqref{eq:coucou}, we know that $\nabla_{v} P$ is $\delta$ smoothing but $u\mapsto  2 (\partial_\mu P(u) ) u $ is clearly not smoothing. Nevertheless, since $\partial_\mu P(u) \in \mathbb{R}$, it generates a trivial flow (a gauge transform) which commutes with everything. More precisely, if $u$ is a solution to \eqref{eq:NLS}, 
$$
w := \exp\Big( 2 \ic \int_0^t \partial_\mu P(u(\mathfrak{t})) \mathrm{d}\mathfrak{t} \Big)\, u \quad \mathrm{and} \quad z := \tau(w),
$$
then $z$ is solution to 
$$
\ic \partial_t z + \partial_x^2 z = \nabla_v P (z).
$$
The nonlinearity of this last equation is smoothing but it is no more Hamiltonian\footnote{note that it is still reversible.}. It means that we have made the nonlinearity of \eqref{eq:NLS} smoothing up to a gauge transform.

\subsection{Remedy} Finally, we explain how we overcome the obstructions described in Subsection \ref{sub:Birk} to get the regularizing normal described in Subsection \ref{sub:desc_reg}. 

\subsubsection*{Convergence of the normal form procedure} To make the normal form convergent, it suffices to  remove only the terms with large divisors in the Birkhoff normal form procedure. More precisely, if we remove only the terms $u_{\boldsymbol{\ell}_1}^{\boldsymbol{\sigma}_1} \cdots  u_{\boldsymbol{\ell}_{2q}}^{\boldsymbol{\sigma}_{2q}}$ such that
$$
 |\boldsymbol{\sigma}_1 \boldsymbol{\ell}_1^2+ \cdots +\boldsymbol{\sigma}_{2q}  \boldsymbol{\ell}_{2q}^2  | \gtrsim q^9
 $$
 then it is not too difficult to prove that the normal form procedure become convergent (the exponent $9$ is not optimal). In other words, there exists a symplectic change of variables $\tau$ close to the identity such that
 $$
 H \circ \tau^{-1} = Z_2+B
 $$
 where $B=G+M$ admits the same decomposition as before except that
 \begin{itemize}
 \item now $r=+\infty$, so $B$ is no more a polynomial but an analytic function,
 \item the condition $ \boldsymbol{\sigma}_1 \boldsymbol{\ell}_1^2+ \cdots +\boldsymbol{\sigma}_{2q}  \boldsymbol{\ell}_{2q}^2 = 0$ has to be replaced by $|\boldsymbol{\sigma}_1 \boldsymbol{\ell}_1^2+ \cdots +\boldsymbol{\sigma}_{2q}  \boldsymbol{\ell}_{2q}^2 | \lesssim q^9$.
 \end{itemize}
 Nevertheless, \eqref{eq:ca_implique} is still valid in this case and so $G$ is still as smoothing as before.
 
\subsubsection*{Terms in $M$ coming from $H$} Now to prove that $B$ is smoothing it suffices to prove that $M$ is smoothing in some sense. Of course $M$ still contains the terms \eqref{eq:relou} and so cannot be smoothing. So, as we discussed in the previous subsection, from now we allow the coefficients of the Hamiltonian to depend on the $L^2$ norm. It completely solves the problem we have for the integrable terms. Indeed, basic results of commutative algebra\footnote{i.e. the fact that the integrable terms are symmetric polynomials of the actions.} allow to prove that the integrable terms of $H$ belong to the algebra generated by the Newton sums
$$
\sum_{k\in \mathbb{Z}} |u_k|^{2p}, \quad p \in \mathbb{Z}.
$$
Thus, it suffices to note that among these Newton sums, the $L^2$ norm is the only non-smoothing one.
Nevertheless, there are terms in $M$ coming from $H$ which are not integrable. 

\medskip

This is where Wick renormalization comes into play. It provides the composition
\begin{equation}
\label{eq:ca_peut_etre_sympa_a_citer}
\| u \|_{L^{2p}}^{2p} = \sum_{\substack{0\leq q\leq p\\ q\neq 1}} \| u\|_{L^2}^{2(p-q)} \frac{p!}{q!}{\binom p q}W_{2q}(u) 
\end{equation}
where
$$
W_{2q}(u) := \int_{\mathbb{T}}\,:\!|u|^{2p}\!:\,(x)\frac{\mathrm{d}x}{2\pi} \quad \mathrm{with} \quad :\!|u|^{2p}\!:\,=\sum_{q=0}^p(-1)^{p-q}\frac{p!}{ q!}{\binom p q}\|u\|_{L^2}^{2(p-q)}|u|^{2q}.
$$
The point is that, rewriting $W_{2q}$ in the Fourier basis as
$$
W_{2q}(u)  = \sum_{\substack{\boldsymbol{\sigma}_1 + \cdots +\boldsymbol{\sigma}_{2q} = 0 \\ \boldsymbol{\sigma}_1 \boldsymbol{\ell}_1+ \cdots +\boldsymbol{\sigma}_{2q}  \boldsymbol{\ell}_{2q} =0 }  } W_{2q}^{\boldsymbol{\ell},\boldsymbol{\sigma}} u_{\boldsymbol{\ell}_1}^{\boldsymbol{\sigma}_1} \cdots  u_{\boldsymbol{\ell}_{2q}}^{\boldsymbol{\sigma}_{2q}} \quad \mathrm{with} \quad W_{2q}^{\boldsymbol{\ell},\boldsymbol{\sigma}}\in \mathbb{R}, \ |W_{2q}^{\boldsymbol{\ell},\boldsymbol{\sigma}}| \leq 2^q q!,
$$
it satisfies the property that any paring is over-paired: if there exist $1\leq i<j\leq 2q$ such that $\boldsymbol{\ell}_i=\boldsymbol{\ell}_j$ and $\boldsymbol{\sigma}_i\boldsymbol{\sigma}_j=-1$ then there exits $1\leq k\leq 2p$ with $k\notin\{i,j\}$ such that $\boldsymbol{\ell}_k=\boldsymbol{\ell}_i=\boldsymbol{\ell}_j$. In particular, all the terms in $W_{2q}$ whose highest index is associated with an action are smoothing. Thus, since we allow the coefficients to depend on the $L^2$ norm, as a consequence of \eqref{eq:ca_peut_etre_sympa_a_citer}, the terms of $M$ coming from $H$ are smoothing.

\subsubsection*{Terms in $M$ generated by the normal form procedure} The discussion is not over yet, as terms in $M$ can be generated by the change of variable. For example, when putting the cubic NLS in Birkhoff normal form, we generate the sixth order integrable term
$$
Z_6(u) = \sum_{k\neq \ell} \frac{|u_k|^2 |u_\ell|^4}{(k-\ell)^2}.
$$
We refer for example to \cite{BFG20} for a proof. Note that this Hamiltonian is only $2$-smoothing: its regularizing effects do not depend on the smoothness of $u$. 

\medskip

So we have to understand what's happen when a new term with a high action is generated by the normal form procedure.  We only remove terms whose  two largest modes are not paired (else we know yet that they are smoothing). They are of the form
\begin{equation}
\label{eq:cond_sympa}
u_{\boldsymbol{\ell}_1}^{\boldsymbol{\sigma}_1} \cdots  u_{\boldsymbol{\ell}_{2q}}^{\boldsymbol{\sigma}_{2q}}, \quad \mathrm{with} \quad (\boldsymbol{\ell},\boldsymbol{\sigma}) \notin \mathbf{P}_q.
\end{equation}
The divisor associated with this term is
$$
\Omega_{\boldsymbol{\ell},\boldsymbol{\sigma}} = \sum_{1\leq j\leq 2q} \boldsymbol{\sigma}_j \boldsymbol{\ell}_j^2.
$$
It is a classical exercise to check that, since $\boldsymbol{\ell},\boldsymbol{\sigma}$ satisfies the zero momentum condition and $(\boldsymbol{\ell},\boldsymbol{\sigma}) \notin \mathbf{P}_q$, we have
\begin{equation}
\label{eq:smoothing_effect}
|\Omega_{\boldsymbol{\ell},\boldsymbol{\sigma}} |\gtrsim_q \frac{|\boldsymbol{\ell}_1^*|}{|\boldsymbol{\ell}_3^*|^2}.
\end{equation}
Note that it implies that the vector field associated with the change of variable is $1$-smoothing (up to a gauge transform). Then it suffices to understand why monomials  $ u_{\boldsymbol{h}_1}^{\boldsymbol{\varphi}_1} \cdots  u_{\boldsymbol{h}_{2q''}}^{\boldsymbol{\varphi}_{2q''}}$ with $(\boldsymbol{h},\boldsymbol{\varphi}) \in \mathbf{P}_{q''}$ generated by Poisson brackets of the form\footnote{where both monomials satisfy implicitly the zero momentum condition.}
$$
\{ u_{\boldsymbol{\ell}_1}^{\boldsymbol{\sigma}_1} \cdots  u_{\boldsymbol{\ell}_{2q}}^{\boldsymbol{\sigma}_{2q}} , u_{\boldsymbol{k}_1}^{\boldsymbol{\varsigma}_1} \cdots  u_{\boldsymbol{k}_{2q'}}^{\boldsymbol{\varsigma}_{2q'}}\}
$$
 satisfy
$$
|\Omega_{\boldsymbol{\ell},\boldsymbol{\sigma}}|\gtrsim_q \frac{|\boldsymbol{h}_1^*|}{  |\boldsymbol{h}_3^*|^2}. 
$$
It is not entirely obvious, but it is not very difficult either (one has to properly use the zero momentum condition and to consider different cases). In any case, this proves that the new terms generated by the change of variable, whose highest mode is in the form of action, are 1-smoothing (up to a gauge transform).

\section{ The KAM part or how to deal with internal parameters}\label{sec:scheme}
In this section we give a scheme of the proof of our result avoiding technicalities. We will use the regularization normal form described in the previous section as a black box and we focus on the KAM procedure.
Without loss of generality, we assume here that $f'(0)=1$. 

\medskip

 We want to explain how we implement the P\"oschel scheme using internal parameters (which are essentially the squared moduli of the initial datum's Fourier coefficients in the final coordinates).  In a first step (see section \ref{sec:step2} and section \ref{sec:step3} above) we construct   invariants tori of  finite dimension in essentially the same way than Kuksin-P\"oschel have done in \cite{KP96}.
Then we implement a procedure that consists in adding one site to the finite-dimensional tori already constructed and we iterate. This is precisely where we need to be very meticulous if we want to take the limit to obtain an infinite-dimensional invariant torus.
By the way we prove more than Theorem \ref{thm:0}: we construct a large family of infinite dimensional tori that  accumulate on  finite dimensional invariant tori (those built  by Kuksin--P\"oschel in \cite{KP96}).

  \medskip

\subsection{Step $1$: set-up.}
Our regularization normal form result provides a symplectic  change of variable $\tau$ defined on a neighborhood of the origin in $H^s$ and close to the identity such that
$$
H^{\mathrm{(reg)}}(u)=H^{\eqref{eq:NLS}} \circ \tau (u)=  H^{(0)}_{1}(u)   +  \underbrace{\sum_{2n+q \geq 6} \sum_{\substack{\boldsymbol{\sigma} \in \{-1,1\}^q \\ \boldsymbol{\sigma}_1+ \cdots + \boldsymbol{\sigma}_q = 0}} \sum_{ \boldsymbol{\sigma}_1 \boldsymbol{\ell}_1+\cdots+  \boldsymbol{\sigma}_q \boldsymbol{\ell}_q =0} H_{n}^{\boldsymbol{\ell},\boldsymbol{\sigma}} \|u\|_{L^2}^{2n} \prod_{1\leq j \leq q} u_{\boldsymbol{\ell}_j}^{\boldsymbol{\sigma}_j}}_{=:R}
$$
with the usual convention $ u_{\ell}^{-1} :=  \overline{u_{\ell}}$ and where
$$
 H^{(0)}_{1}(u)  := \frac12 \sum_{k\in \mathbb{Z}} \big(  |k|^2 -  \frac{|u_k|^2}2 \big) |u_k|^2 + \frac12\| u\|_{L^2}^4
$$
and the coefficients $ H_{n}^{\boldsymbol{\ell},\boldsymbol{\sigma}} \in \mathbb{C}$ satisfy the bound
$$
| H_{n}^{\boldsymbol{\ell},\boldsymbol{\sigma}}| \lesssim C^{q+2n} \big( 1 \wedge \frac{  \langle \boldsymbol{\ell}_3^* \rangle^2  }{  \langle \boldsymbol{\ell}_1^* \rangle } \big),\quad \text{for some constant } C>0
$$
and  where  $\boldsymbol{\ell}^*$ denotes a rearrangement of $\boldsymbol{\ell}$ such that $|\boldsymbol{\ell}^*_1|\geq \cdots \geq |\boldsymbol{\ell}^*_{q}|$.

This bound implies that $H^{\mathrm{(reg)}}$ is 1-smoothing in some week sense discuted in section \ref{sec:reg}. Notice that this expansion depends on $ \|u\|_{L^2}^2$ which is a constant of motion for \eqref{eq:NLS}. Nevertheless, $ \|u\|_{L^2}^2$ will generate a dependency  on our internal parameters (denoted by $\xi$ in the sequel, see \eqref{xi}), so, in principle, we must take this dependence into account (in particular to control Lipschitz norms). It turns out that this dependency does not seriously affect the proof scheme, but it does lead to a clear increase in the complexity of the notations.  So, to simplify our presentation, in this sketch of proof we'll consider $ \|u\|_{L^2}^2$ as a constant.

\medskip

We are working in the vicinity of the origin, so to take advantage of the smallness of the solution we make a change of scale $u\leadsto \eps u$ and we still consider $u$ to be small (if we prefer, we cut this smallness in half). This subterfuge enables us to isolate a first small parameter $\eps$ which will be useful in several places in the proof (but of course a finite number of times). So we set
\begin{align*}
H^{\mathrm{(reg)}}_\eps(u)&:=\eps^{-4}H^{\mathrm{(reg)}}(\eps u)\\
&= H^{(0)}_{\eps}(u) +\eps R_\eps
\end{align*}
where 
$$
 H^{(0)}_{\eps}(u)  := \frac12 \sum_{k\in \mathbb{Z}} \big( \eps^{-2} |k|^2 -  \frac{|u_k|^2}2 \big) |u_k|^2
$$
and, $R$ being of order 6, $R_\eps=O(\eps)$. 

\medskip

\subsection{Step $2$: first opening.}\label{sec:step2}
After this initial set-up, it's time for the first opening. We fix $\S_1$ a finite subset of $\Z$ and we write $u_k,\ k\in \S_1$ in action-variables $(y_k,\theta_k)$, i.e. we set
\begin{equation}\label{xi}
u_k=(\xi_k+y_k)^{\frac12}e^{i\theta_k}, \quad k\in\S_1
\end{equation}
where $(\xi_k)_{k\in\S_1}$ are our first parameters to play with. 
Let us denote
$$\psi_1: (u_k)_{k\in\Z}\mapsto ((y_k)_{k\in\S_1},(\theta_k)_{k\in\S_1};(u_k)_{k\in\S_1^c})$$
 the corresponding change of variable. 
For the moment we assume $r_1^{2\nu}<|\xi|_{\ell^\infty}<2 r_1^{2\nu}$ for some small $r_1$ to be defined and some fixed $1/2<\nu<1$. We work in a neighborhood of the torus
$\mathrm T_\xi$ (see \eqref{eq:def_T_droit}): $$|y_k|<r_1^2,\ k\in \S_1 \text{ and } \sum_{k\in\S_1^c}\varpi_k^2 |u_k|^2<r^2_1.$$
In this new variables the quadratic part of $ H^{(0)}_{\eps}\circ \psi_1$ reads
\begin{equation}\label{Z}
Z^{(0)}_{\eps}=\frac12 \sum_{k\in \S_1} \big( \eps^{-2} |k|^2 +2\xi_k \big)y_k + \frac12 \sum_{k\in \S_1^c}  \eps^{-2} |k|^2  |u_k|^2
\end{equation}
and the quartic part of $H^{\mathrm{(reg)}}_\eps\circ \psi_1$ reads 
$$
Q^{(0)}_{\eps}=\frac12 \sum_{k\in \S_1} y_k^2+ \frac12 \sum_{k\in \S_1^c}    |u_k|^4+O(\eps).$$
Let us introduce some vocabulary. A formal Hamiltonian
$$
P(y,\theta;u) :=\sum_{ \boldsymbol{k} \in \mathbb{Z}^{\mathcal{S}_1}  }\sum_{ \boldsymbol{m} \in \mathbb{N}^{\mathcal{S}_1}  } \sum_{q\in \mathbb{N}} \sum_{ \boldsymbol{\ell} \in (\mathcal{S}_1^c)^q  }  \sum_{ \boldsymbol{\sigma} \in  \{-1,1\}^q  } P_{\boldsymbol{k},\boldsymbol{m}}^{\boldsymbol{\ell},\boldsymbol{\sigma}} \,  \underbrace{ \prod_{i=1}^q u_{\boldsymbol{\ell}_i}^{\boldsymbol{\sigma}_i} }_{=:u_{\boldsymbol{\ell}}^{\boldsymbol{\sigma}} } \, \underbrace{ \prod_{j\in \mathcal{S}_1} y^{\boldsymbol{m}_j}_j e^{\ic\theta_j\boldsymbol{k}_j}  }_{=: y^{\boldsymbol{m}} \, e^{ \ic\theta\cdot\boldsymbol{k} }} \, 
$$
is said to be
\begin{itemize}
\item  integrable if $P_{\boldsymbol{k},\boldsymbol{m}}^{\boldsymbol{\ell},\boldsymbol{\sigma}} \neq0$ only when $\boldsymbol{k}=0$ and there exists a permutation $\varphi$ such that $\varphi\boldsymbol{\ell}=\boldsymbol{\ell}$ and $\varphi\boldsymbol{\sigma}=-\boldsymbol{\sigma}$
\item an adapted jet  if $P$ does not contains integrable terms and $P_{\boldsymbol{k},\boldsymbol{m}}^{\boldsymbol{\ell},\boldsymbol{\sigma}} \neq0$  only when $|\boldsymbol{m}|_{\ell^1}\leq2$ and $q\leq3$.
\item in normal form if it is integrable and at most quartic: $P_{\boldsymbol{k},\boldsymbol{m}}^{\boldsymbol{\ell},\boldsymbol{\sigma}} \neq0$ only when $\boldsymbol{k}=0$, there exists a permutation $\varphi$ such that $\varphi\boldsymbol{\ell}=\boldsymbol{\ell}$ and $\varphi\boldsymbol{\sigma}=-\boldsymbol{\sigma}$ and $2|\boldsymbol{m}|_{\ell^1}+q\leq4$.
\end{itemize}
We then denote by $\Pi_{\mathrm{int}}$, $\Pi_{\mathrm{ajet}} $, $\Pi_{\mathrm{nor}}$ the corresponding projections. Finally we introduce $\Pi_{\mathrm{rem}}$ the projection defined by $\Pi_{\mathrm{rem}}P=P- \Pi_{\mathrm{nor}}P-\Pi_{\mathrm{ajet}} P$. We note that these notions are linked to $\S_1$ and  sometimes (when necessary for clarity) we will denote $\Pi^{\S_1}_{\mathrm{ajet}} $.

\medskip

\subsection{Step $3$: first KAM procedure, \`a la Kuksin--P\"oschel.}\label{sec:step3}
At this stage there is nothing really original. Like Kuksin-P\"oschel in \cite{KP96}, we use the term $\xi_k y_k$ in \eqref{Z} to modulate the frequencies of the quadratic part of the Hamiltonian $H^{\mathrm{(reg)}}_\eps$ and the smallness of the 2-jet of the remainder part for $\eps$ small enough  to obtain for $\xi$ in a set $O_\eps^{(1)}$ of asymptotically full measure when $\eps\to 0$, that, for $\eps$ small enough, the linear torus $\mathrm T_\xi$ persists, slightly deformed, as invariant sets of \eqref{eq:NLS} . \\
But we want to be able to work around these new tori and thus we need more. 
In fact we prove that, for $\xi\in O_\eps^{(1)}$, there exists a symplectic change of variable $\tau_1$ close to identity such that 
$$H^{\mathrm{(reg)}}_\eps\circ\psi_1\circ\tau_1= H^{(0)}_{\eps}\circ\psi_1+\eps K_\eps^{(1)}\quad \text{ with } \Pi_{\mathrm{ajet}}K_\eps^{(1)}=0.$$
Furthermore $\|\Pi_{\mathrm{rem}}K_\eps^{(1)}\|\leq 1$ and $\| \Pi_{\mathrm{nor}}(K_\eps^{(1)}-R_\eps\circ\psi_1)\|\leq r_1$. In particular in the new variables the frequencies read
\begin{equation}\label{om}\omega_j^{(1)}=\eps^{-2}j^2+2 \mathbbm{1}_{j\in \mathcal{S}_1} \xi_j+\eps\lambda_j^{(1)}(\xi),\quad j\in \Z
\end{equation}
where $\lambda_j^{(1)}$ is 1-Lipschitz and bounded by 1.

\medskip

Some comments:
\begin{itemize}
\item  To obtain such normal form, we need to control the so-called small divisors. Concretely we prove: for all $\xi \in O_\eps^{(1)}$, all $\boldsymbol{k}\in \mathbb{Z}^{\mathcal{S}_1}\setminus \{0\}$, all $d \in \llbracket 0,4 \rrbracket$, all vector $\boldsymbol{\ell} \in (\mathcal{S}_1^c)^d$, all $\boldsymbol{b} \in (\mathbb{Z}^*)^d$,  provided that
$\|\boldsymbol{b}\|_{\ell^1}$ not to large, 
\begin{equation}\label{smalldiv}\Big|  \sum_{i\in \mathcal{S}_1} \boldsymbol{k}_i \omega^{(0)}_i(\xi) + \sum_{p=1}^d \boldsymbol{b}_p \omega_{\boldsymbol{\ell}_p}(\xi) \Big| \gtrsim\left(\frac{r_1^{2\nu}}{ \| \boldsymbol{k} \|_{\ell^\infty}^{\# \mathcal{S}_1} }\right)^{\alpha_d}
\end{equation}
for some constant $\alpha_d$. The presence of $r_1^{2\nu}$ in the numerator of the right hand side is due to the size of our parameters $\xi$.
We note that we have 5 estimates depending on the value of $d$. For $d=1$ and $d= 2$ they correspond to the so called first and second Melnikov conditions. By analogy, for $d=3$ and $d= 4$ we can see these conditions as third and fourth Melnikov conditions.  
We prove that they are typically satisfied using arguments similar to those we developed  in \cite{BG22} to perform Birkhoff normal forms for PDEs in low regularity (convexity arguments could also be used). 
\item These third  Melnikov condition is needed  to remove the adapted jet during the KAM step. This is very expensive in terms of small divisors. In fact to perform such KAM theorem, we have to solve successively 6 cohomological equations and at the end we have to assume that this part of the adapted jet of our perturbation is very small ($\lesssim r_1^{4000}$ in \cite{BGR24}). 
\item Since $R_\eps=O(\eps)$, we can get rid of this factor  $r_1^{4000}$ at the denominator by choosing $\eps\lesssim r_1^{5000}$.
\item Only at this price  we can control the normal form  of our new Hamiltonian $\Pi_{\mathrm{nor}}K_\eps^{(1)}$. As we shall see, this will be essential, after a new opening, to obtain new frequencies that still satisfy an estimate of the type \eqref{om} with still a Lipschitz control on $\lambda_j(\xi)$ (see \eqref{om2}).

\end{itemize}

\medskip

\subsection{Step $4$: second opening and Birkhoff procedure.} We want to open a new site (only one at a time, it will be important), say $\S_2=\S_1\cup\{i_2\}$ with $i_2\notin \S_1$. This generates a new  change of variable $\psi_2$ and we set
$$H^{\mathrm{(reg)}}_\eps\circ\psi_1\circ\tau_1\circ\psi_2= H^{(0)}_{\eps}\circ\psi_1\circ\psi_2+\eps K_\eps^{(1)}\circ\psi_2.$$
The new internal frequencies read 
\begin{equation}\tilde\omega_j^{(1)}=\eps^{-2}j^2+2 \mathbbm{1}_{j\in \mathcal{S}_2} \xi_j+\eps\tilde\lambda_j^{(1)}(\xi),\quad j\in \Z
\end{equation}
where $\tilde\lambda_j^{(1)}$ is still Lipschitz and bounded. Assume $r_2^{2\nu}<|\xi_{i_2}|<2 r_2^{2\nu}$ for some small $r_2\leq r_1$ to be defined. We will now work in a neighborhood of the torus\footnote{Where given $\xi\in\ell^1_s$ and $\S\subset\Z$ we denote $\xi^{(\S)}$ the element of $\R^\Z$ defined by $\xi^{(\S)}_i=\xi_i,\ i\in\S$ and $\xi^{(\S)}_i=0,\ i\in\S^c$.}
$\mathrm T_{\xi^{(\S_2)}}$:   
\begin{equation}\label{neighborhood2}|y_k|<r_2^2,\ k\in \S_2 \text{ and } \sum_{k\in\S_2^c}\varpi_k^2 |u_k|^2<r^2_2.\end{equation}
In fact $r_2$  represents the distance  to the $\mathrm T_{\xi^{(\S_1)}}$, the invariant torus constructed in the previous step, and so we want to choose $r_2$ as small as we want independently of $\eps$.

We can still modulate internal frequencies in such a way  estimate \eqref{smalldiv}, with $\S_1$ replaced by $\S_2$ and $r_1$ replaced by $r_2$,  be satisfied in a "large" set $O_\eps^{(2)}$. Nevertheless we cannot use again the smallness of $\eps$ because we want to iterate infinitely many times this procedure of adding a new site, generating a decreasing sequence of radii $r_p,\ p\geq2$, and  we cannot impose constraints of the type $\eps\lesssim r_p^{5000}$ for all $p\geq2$.

 So our small parameter is now $r_2$ and thus  \eqref{smalldiv} with $r_1$ replace by $r_2$ is not a good control at all. The idea is that we can still get  \eqref{smalldiv} with $r_1$ if the small divisor (or if you prefer the corresponding monomials, $y^{\boldsymbol{m}} \, e^{ \ic\theta\cdot\boldsymbol{k} } u_{\boldsymbol{\ell}}^{\boldsymbol{\sigma}}$, $\boldsymbol{k}\in \Z^{\S_2}$, $\boldsymbol{\ell}\in (\S_2^c)^d$) contains one internal site in $\S_1$ and $\boldsymbol{k}_{i_2}$ is not too large, i.e. if one the $k_i$, $i\in\S_1$, does not vanish and $|\boldsymbol{k}_{i_2}|$ is smaller than a given constant. In that case,  we are done for the following reasons:
\begin{itemize}
\item The monomials that we want to eliminate comes, at the previous step, from a term  at least of order 4 in $(y,u)$ (we killed the 3-jet).
 We work on the neighborhood \eqref{neighborhood2} and thus this monomials has a size $\lesssim r_2^4$ and its vectorfield has a size $\lesssim r_2^2$ .
\item In the resolution of the cohomological equation we have to control this vectorfield divided by the small denominator  (see \eqref{smalldiv}), so something like $\lesssim \frac{r_2^2}{r_1^{2\nu\alpha_3}}$.
\item Even if we have to solve successively 6 such equations (as explain before), we will face something like $\lesssim \frac{r_2^2}{r_1^{12\nu\alpha_3}}$.
\item So we just have to choose $r_2$ small enough compare to $r_1$.
\end{itemize}

Now it remains to consider monomials of the kind $y^{\boldsymbol{m}} e^{ \ic\theta\cdot\boldsymbol{k} }  u^{a}_{i_2}\bar u^{\bar a}_{i_2}u_{\boldsymbol{\ell'}}^{\boldsymbol{\sigma'}} $ with $\boldsymbol{k}\in \Z^{\S_1}$, $\boldsymbol{\ell'}\in(\S_2^c)^{d'}$, $\boldsymbol{\sigma'}\in\{-1,1\}^{d'}$ and $a,\bar a\in\N$. Since we are interested in eliminated the adapted jet, we can restrict the study to the case $d'\leq 3$. Indeed, since $|u_{i_2}|\lesssim r_2^\nu$, we can compensate a small divisor of order $r_2^{4000}$ (generated, as in the previous step, by  \eqref{smalldiv} with $r_1$ replaced by $r_2$) if  $a+\bar a\geq 10000$ (recall that $\nu>1/2$). So we say that such monomials is "to be removed" if it is not integrable and $d'\leq3$ and $a+\bar a\leq 10000$. We denote by $\Pi_{\mathrm{tbr}}$ the associated projection.

To make the "to be removed" part of the Hamiltonian much smaller, we apply  a Birkhoff procedure {\it before} the second opening.
To apply a Birkhoff procedure we have to check that we have a good control the small divisors associated to these "to be removed" monomials. It turns out that  they are non resonant and thus the associated small divisor depends again on $r_1$, so not dangerous.  

After the Birkhoff step, i.e. applying a symplectic change of variable $\varphi_1$ and removing a small part of parameters (we will ignore it in this scheme), we get the existence of $r_2\ll r_1$ such that on a neighborhood of $\mathrm T_{\xi^{(\S_2)}}$ of size $r_2$ (see \eqref{neighborhood2})
$$H^{\mathrm{(reg)}}_\eps\circ\psi_1\circ\tau_1\circ\varphi_1=H^{(0)}_{\eps}\circ\psi_1+\eps K_\eps^{(1)}\circ\varphi_1\quad \text{with } \|\Pi_{\mathrm{tbr}}  K_\eps^{(1)}\circ\varphi_1\|\leq r_2^{5000}$$
in such a way after the second opening we have 
$$H^{\mathrm{(reg)}}_\eps\circ\psi_1\circ\tau_1\circ\varphi_1\circ\psi_2= H^{(0)}_{\eps}\circ\psi_1\circ\psi_2+\eps \tilde K_\eps^{(1)} $$
with $\tilde K_\eps^{(1)}$ having a very small adapted jet, says  $\|\Pi_{\mathrm{ajet}}^{\S_2}  \tilde K_\eps^{(1)}\|\leq r_2^{4000}$, in such a way we can apply a KAM procedure but this time with $r_2$ playing the role of small parameter.

\subsection{Step $5$ : second KAM procedure.} With a maybe smaller $r_2$, we prove the existence of a good set $O_\eps^{(2)}$ with $\frac{\mathrm{meas}(O_\eps^{(1)}\setminus O_\eps^{(2)})}{\mathrm{meas}(O_\eps^{(1)})}=O(r_2^{1/1000})$ and 
there exists a symplectic change of variable $\tau_1$ close to identity such that for $\xi\in O_\eps^{(2)}$
$$H^{\mathrm{(reg)}}_\eps\circ\psi_1\circ\tau_1\circ\varphi_1\circ\psi_2\circ\tau_2= H^{(0)}_{\eps}\circ\psi_1+\eps K_\eps^{(2)}\quad \text{ with } \Pi_{\mathrm{ajet}}^{\S_2}K_\eps^{(2)}=0.$$
Furthermore $\|\Pi_{\mathrm{rem}}^{\S_2}K_\eps^{(2)}\|\leq 1$ and $\| \Pi^{\S_2}_{\mathrm{nor}}(K_\eps^{(2)}-K^{(1)}_\eps\circ\psi_2)\|\leq r_2$. In particular the new internal frequencies read
\begin{equation}\label{om2}\omega_j^{(2)}=\eps^{-2}j^2+2 \mathbbm{1}_{j\in \mathcal{S}_2} \xi_j+\eps\lambda_j^{(2)}(\xi),\quad j\in \Z
\end{equation}
where $\lambda_j^{(2)}$ is still Lipschitz and bounded.

\subsection{Step $6$ : iteration of the loop.}
Now we can iterate this loop, Birkhoff-opening-KAM, to generate tori of arbitrary dimension $p$ close to $\mathrm T_{\xi^{(\S_1)}}$. Of course it remains to check that this scheme converges when $p\to \infty$. Very schematically, if we denote by $(r_p)_{p\in\N}$ the decreasing sequence of the sizes of the circles that make up the infinite dimensional tori,  we strongly use, as in \cite{Pos02}, that $r_{p+1}$ can be chosen much smaller than $r_p$. In a way we have an infinite reservoir of smallness and this enables us to make assumptions such as "the adapted jet is of size $r_p^{4000}$" in our KAM theorem. The downside is that we relinquish all reasonable control over the decay of $r_p$ (contrary to \cite{BMP20} or \cite{Con23}).

\subsection{Additional comments}
$ $

 \noindent \underline{$\triangleright$ Modulation of frequencies and twist condition.} The crucial hypothesis in any KAM theorem is the twist condition: assume that the frequencies vector $\omega=(\omega_j)_{j\in\mathcal S}$, $\mathcal S$ being finite or infinite subset of $\Z$, depends smoothly on parameters $\xi\in\R^{\mathcal S}$, we want to be able to modulate the frequencies when moving the parameters. So we want the derivative of $\omega$, denoted $\mathrm{d}\omega$, to be  invertible in some sense depending on the topology we choose on the set of frequencies and on the set of parameters. 
Following the above scheme of proof, let us write $\omega^{(p)}_j=\xi_j+\eps\lambda^{(p)}_j,\ j\in\S_p$. The twist condition, after applying $p$ loops as described above, results in invertibility of $\text{Id}+\eps \mathrm{d}\lambda^{(p)}$ and, if we want a diagonally dominant matrix\footnote{which is equivalent to consider  $\omega:\ell^1(\mathcal S_p)\to\ell^1(\mathcal S_p)$ and to apply a Neumann series argument}, 
can be rewritten as
\begin{equation}
\label{eq:bound_l1}
\sup_{k\in\mathcal S_p}\sum_{j\in\mathcal S_p} |\partial_{\xi_k}\lambda^{(p)}_j |\lesssim 1.
\end{equation}

The non-linearity regularization improves the estimate on $\lambda^{(p)}$.  By $\delta$-regularizing Hamiltonian we mean an analytic Hamiltonian $P$ that satisfies (locally) $\nabla P: \ \ell^1_{s_0}\to \ell^1_{s_0+\delta}$ where 
$$
\ell^1_{s_0}:=\{ u \in \mathbb{C}^{\mathbb{Z}} \mid \sum_{k\in \mathbb{Z}} \langle k\rangle^{s_0}|u_k|<\infty  \}, \quad \mathrm{for\ some\ } s_0\geq0.
$$ 
Without regularization ($\delta=0$), we  expect estimates on $\lambda^{(p)}$ of the type $|\partial_{\xi_k}\lambda^{(p)}_j |\lesssim 1$ and we cannot expect to prove the twist condition when $p\to\infty$. Thus without regularization we obtain the result of \cite{KP96}.
If we want to be able to consider the limit $p\to\infty$, we need regularization ($\delta>0$) and we can still distinguish two cases:
\begin{itemize}
\item with a regularization without loss ($\delta>0$ and  $s_0=0$) we  expect estimates on $\lambda$ of the type $|\partial_{\xi_k}\lambda^{(p)}_j |\lesssim \frac1{\langle j\rangle^\delta}$ uniformly in $k$ (and $p$) and we can  prove the twist  condition when $p\to\infty$ for $\delta>1$ or $S$ sparse.
\item with a regularization with loss ($\delta>0$ and  $s_0>0$), which is actually our case with $s_0=2\delta=2$, we prove estimates on $\lambda^{(p)}$ of the type\footnote{in fact $\xi\mapsto\lambda$ will be only Lipschitz and we have to adapt the twist condition but we omit this problem here}
\begin{equation}
\label{estim:omeg}
\left|\partial_{\xi_k}\lambda^{(p)}_j\right|\lesssim  (1\wedge \langle j \rangle^{-\delta} \langle k\rangle^{2s_0} \wedge  \langle k \rangle^{-\delta} \langle  j \rangle^{2s_0}  ) \vee  \langle j\rangle^{-\delta}
\end{equation}
 and we prove the twist condition for $\delta>0$ and $p\to\infty$ provided that $\S_\infty=\lim_{p\to\infty} \S_p$ sufficiently sparse. Note that putting together \eqref{eq:bound_l1} and \eqref{estim:omeg},  with $s_0=2\delta=2$, we get the sparsity condition 
 \begin{equation}
\label{eq:def_sparcity_condition_intro}
\sup_{k\in \mathcal{S}_\infty} \sum_{j\in \mathcal{S}_\infty} \Big( 1\wedge \frac{\langle  k \rangle^{4}}{\langle j \rangle}  \wedge  \frac{\langle  j \rangle^{4}}{\langle k \rangle}   \Big) \vee  \frac1{\langle j\rangle} < \infty.
\end{equation}

\end{itemize}

\medskip

 \noindent \underline{$\triangleright$ About the sparsity condition \eqref{eq:def_sparcity_condition_intro} on $\mathcal S_\infty$.} The explicit condition \eqref{eq:def_sparcity_condition_intro} is the only one we impose on the sites.  We can try to explain why estimate \eqref{estim:omeg} (implying \eqref{eq:def_sparcity_condition_intro}) is natural and thus  why, for the meantime, sparsity seems unavoidable without a significant improvement in our regularizing normal form theorem. In fact the problem comes from integrable terms in the Hamiltonian.  Our regularization  normal form allows terms of the form $c_{k,\ell,j}|u_j|^2|u_\ell|^2 |u_k|^2$ in the nonlinear term with, if we assume $\langle k\rangle \gg \langle\ell\rangle\geq \langle j\rangle$,  the estimate\footnote{Indeed we can compute explicitly all the sixtic terms and we obtain (see for instance \cite{BFG20} where these terms are very useful) $Z_6=\sum_{k\neq\ell}\frac{|u_k|^4 |u_\ell|^2}{(k-\ell)^2}$ which leads to a quite better estimates of $c_{k,\ell,j}$, but it is difficult to expect such "convolution" structure to be stable at higher degree.  } $|c_{k,\ell,j}|\leq \langle k\rangle^{-\delta}\langle \ell\rangle^{4\delta}\langle j\rangle^{4\delta}$. After opening modes $j$ and $\ell$ this term generates a contribution to the frequency of index $k$: $\tilde\omega_k(\xi)=c_{k,\ell,j}\xi_\ell\xi_j$ and we have $\partial_{\xi_\ell} \tilde\omega_k =c_{k,\ell,j}\xi_j$. If $j$ is very small, this term could be of order $\langle k\rangle^{-\delta}\langle \ell\rangle^{4\delta}$.

%

\medskip
 
 \noindent \underline{$\triangleright$ About the amplitudes of the tori.} Very schematically, if we denote by $(r_p)_{p\in\N}$ the sequence of sizes of the circles that make up the infinite dimensional tori we are constructing, we strongly use, as in \cite{Pos02}, that $r_{p+1}$ can be chosen much smaller than $r_p$. In a way we have an infinite reservoir of smallness and this enables us to make assumptions such as "the adapted jet is of size $r^{4000}$" in our KAM theorem. The downside is that we relinquish all reasonable control over the decay of $r_p$ (contrary to \cite{BMP20} or \cite{Con23}). More generally, we are not trying to optimize sizes.


\begin{thebibliography}{100000000}




\bibitem[BBHM18]{BBHM18}
{\sc P.~Baldi, M.~Berti, E.~Haus, R.~Montalto}, 
\emph{Time quasi-periodic gravity water waves in finite depth}, \href{https://doi.org/10.1007/s00222-018-0812-2}{Inventiones Math}.
\textbf{214} (2018), 739--911.





%
\bibitem[BC24]{BC24}
{\sc J.~Bernier, N.~Camps},  
\emph{Long time stability for cubic nonlinear Schr\" odinger equations on non-rectangular flat tori},
  \href{https://arxiv.org/abs/2402.04122}{arXiv:2402.04122}, (2024)
%
%
%
%
 \bibitem[BFG20]{BFG20}
{\sc J.~Bernier, E.~Faou, B.~Gr\'ebert}, 
\emph{Rational normal forms and stability of small solutions to nonlinear Schr\"odinger equations}, 
\href{https://doi.org/10.1007/s40818-020-00089-5}{Annals of PDE} {\bf 6}, 14 (2020).
%
%
%
 \bibitem[BG21]{BG21}
{\sc J.~Bernier, B.~Gr\'ebert}, 
\emph{Long time dynamics for generalized Korteweg-de Vries and Benjamin-Ono equations,} 
\href{https://doi.org/10.1007/s00205-021-01666-z}{ Arch. Rational Mech. Anal.} \textbf{241}, 1139--1241 (2021)

\bibitem[BG22]{BG22}
{\sc J.~Bernier, B.~Gr\'ebert}, 
\emph{Birkhoff normal forms for Hamiltonian PDEs in their energy space}, 
\href{https://doi.org/10.5802/jep.193}{Journal de l’Ecole polytechnique — Mathématiques}, Tome 9 (2022), pp. 681-745.
%
%
%
%

\bibitem[BGR24]{BGR24}
{\sc J.~Bernier, B.~Gr{\'e}bert, T.~Robert},  
\emph{Infinite dimensional invariant tori for nonlinear Schr\"odinger equations},
  \href{https://arxiv.org/abs/2412.11845}{arxiv.org/abs/2412.11845}, (2024)





\bibitem[BB13]{BB13}
{\sc M.~Berti, P.~Bolle},
\emph{Quasi-periodic solutions with Sobolev regularity of NLS on $\mathbb T^d$ with a multiplicative potential},
\href{https://doi.org/10.4171/JEMS/361}{J. Eur. Math. Soc.} \textbf{15} (2013), 229--286.



\bibitem[BHM23]{BHM23} 
{\sc M.~Berti, Z.~Hassainia, N.~Masmoudi},
\emph{Time quasi-periodic vortex patches of Euler equation in the plane},
\href{https://doi.org/10.1007/s00222-023-01195-4}{Inventiones Math.}, \textbf{233} (2023), 1279--1391.
	


\bibitem[BKM18]{BKM18}
{\sc M.~Berti, T.~Kappeler, R.~Montalto},
\emph{Large KAM tori for perturbations of the dNLS equation},
\href{https://doi.org/10.24033/ast.1053}{Asterisque}, 403. 2018.	


\bibitem[BCGP24]{BCGP24}
{\sc L.~Biasco, L.~Corsi, G.~Gentile, M.~Procesi},
\emph{Asymptically full measure sets of almost-periodic solutions for the NLS equation},
\href{https://arxiv.org/abs/2412.02648}{arxiv:2412.02648}

 
\bibitem[BMP20]{BMP20}
{\sc L.~Biasco, J.E.~Massetti, M.~Procesi},
\emph{Almost periodic invariant tori for the NLS on the circle},
\href{https://doi.org/10.1016/j.anihpc.2020.09.003}{Ann. Inst. H. Poincar\'e C Anal. Non Lin\'eaire}, \textbf{38} (2020) no.~3, 711--758.


\bibitem[BMP23]{BMP24}
{\sc L.~Biasco, J.E.~Massetti, M.~Procesi},
\emph{Weak Sobolev almost periodic solutions for the 1D NLS},
 \href{https://doi.org/10.1215/00127094-2022-0089}{Duke Math. J.}, 172 (14) 2643 - 2714, 2023. 



\bibitem[Boh47]{Boh47}{
{\sc H.~Bohr},
\emph{Almost Periodic Functions},
Chelsea Publishing Co., New York, 1947.
}



\bibitem[Bou94]{Bou94}
{\sc J.~Bourgain},
\emph{Construction of quasi-periodic solutions for Hamiltonian perturbations of linear equations and applications to nonlinear PDE},
\href{https://doi.org/10.1155/S1073792894000516}{Internat. Math. Res. Notice} \textbf{1994} (1994), no.~11, 475--497.
%
%

\bibitem[Bou96b]{Bou96}
{\sc J.~Bourgain},
\emph{Construction of approximative and almost-periodic solutions of perturbed linear Schr\"odinger and wave equations, }
\href{https://doi.org/10.1007/BF02247885}{Geometric and Functional Analysis} \textbf{6} (1996) 201--230. 
 

\bibitem[Bou98]{Bou98}
{\sc J.~Bourgain},
\emph{Quasi-periodic solutions of {H}amiltonian perturbations of 2{D} linear {S}chr\"odinger equations},
\href{https://doi.org/10.2307/121001 }{Ann. of Math.} (2), \textbf{148} (1998), no.~2, 363--439.

%
\bibitem[Bou00]{Bou00}
{\sc J.~Bourgain},
\emph{On diffusion in high-dimensional Hamiltonian systems and PDE},
\href{https://doi.org/10.1007/BF02791532}{J. Anal. Math.} \textbf{80}, 1--35  (2000) 
%

\bibitem[Bou05]{Bou05}
{\sc J.~Bourgain},
\emph{On invariant tori of full dimension for 1{D} periodic NLS},
\href{https://doi.org/10.1016/j.jfa.2004.10.019}{J. Funct. Anal.} \textbf{229} (2005), no.~1, 62--94.

%
%




\bibitem[Con23]{Con23}
{\sc H.~Cong},
\emph{ The existence of full dimensional KAM tori for nonlinear Schr\"odinger equation},
 \href{https://doi.org/10.1007/s00208-023-02782-9}{Math. Ann.}  Volume 390, pages 671–719, (2024).



%
%

\bibitem[CY21]{CY21}
{\sc H.~Cong, X.~Yuan},
\emph{The existence of full dimensional invariant tori for 1-dimensional nonlinear wave equation},
\href{https://doi.org/10.1016/j.anihpc.2020.09.006}{Ann. Inst. H. Poincar\'e C Anal. Non Lin\'eaire} \textbf{38} (2021), no.3, 759--786.


\bibitem[CGP23]{CGP24}
{\sc L.~Corsi, G.~Gentile, M.~Procesi},
\emph{Almost-periodic solutions to the NLS equation with smooth convolution potentials},
  \href{https://arxiv.org/abs/2309.14276}{arXiv:2309.14276}, (2023)



\bibitem[CW93]{CW93} 
{\sc W.~Craig, C.E.~Wayne}, 
\textit{Newton's method and periodic solutions of nonlinear wave equation}, 
\href{https://doi.org/10.1002/cpa.3160461102}{Comm. Pure  Appl. Math.} \textbf{46} (1993), no.~11, 1409--1498.

%

%


\bibitem[EK10]{EK10}
{\sc L.H.~Eliasson, S.B.~Kuksin},
\emph{KAM for the nonlinear Schr\"odinger equation},
\href{https://doi.org/10.4007/annals.2010.172.371}{Ann. Math.} (2) \textbf{172} (2010), 371--435.


\bibitem[EGK16]{EGK16}
{\sc L.H.~Eliasson, B.~Gr\'ebert,  S.B.~Kuksin},
\emph{KAM for non-linear beam equation},
\href{https://doi.org/10.1007/s00039-016-0390-7}{Geom. Funct. Anal.} \textbf{26} (2016) 1588--1715

\bibitem[ErTz13a]{ErTz1}
{\sc M.B.~Erdo\u gan, N.~Tzirakis,} 
\emph{Talbot effect for the cubic nonlinear Schr\"odinger equation on the torus,} 
\href{https://dx.doi.org/10.4310/MRL.2013.v20.n6.a7}{Math. Res. Lett.}, 20 (2013), 1081–1090.
%
 
\bibitem[ErTz13c]{ErTz3}
{\sc M.B.~Erdo\u gan, N.~Tzirakis}
\emph{Smoothing and global attractors for the Zakharov system on the torus,} 
\href{https://doi.org/10.2140/apde.2013.6.723}{Anal. PDE}, 6 (2013), pp. 723--750



\bibitem[FG24]{FG24}
{\sc R.~Feola, F.~Giuliani},
\emph{Quasi-periodic traveling waves on an infinitely deep perfect fluid under gravity}, 
\href{https://doi.org/10.1090/memo/1471}{Memoirs of the American Mathematical Society}, 158 p. (2024). 



%
%
%
%

%
%
\bibitem[GX13]{GX13}
{\sc J.~Geng, X.~Xu},
\emph{Almost-periodic solutions of one dimensional Schr\"odinger equation with the external
parameters}, \href{https://doi.org/10.1007/s10884-013-9302-9}{J. Dyn. Differ. Equations} \textbf{25} (2013), no. 2, 435--450.



%

%


\bibitem[GP16]{GP16}
{\sc B.~Gr\'ebert, E. Paturel},
\emph{KAM for the Klein Gordon equation on $\mathbb S^d$},
\href{https://doi.org/10.1007/s40574-016-0072-2}{Boll. Unione Mat. Ital}, \textbf{9} (2016), 237--288.
%

\bibitem[KP03]{KP03}
{\sc T.~Kappeler, J.~P{\"o}schel},
\emph{KdV\&KAM}, 
\href{https://doi.org/10.1007/978-3-662-08054-2}{Springer} (2003).

%


\bibitem[Kuk87]{Kuk87}
{\sc S.~Kuksin}, 
\emph{Hamiltonian perturbations of infinite-dimensional linear systems with imaginary spectrum},
\href{https://doi.org/10.1007/BF02577134}{Funktsional. Anal. i Prilozhen.} \textbf{21} (1987), no.~3, 22--37.


%

\bibitem[KP96]{KP96}
{\sc S.~Kuksin, J.~P{\"o}schel},
\emph{Invariant {C}antor manifolds of quasi-periodic oscillations for a nonlinear {S}chr{\"o}dinger equation},
\href{https://doi.org/10.2307/2118656}{Ann. of Math.} (2) \textbf{143} (1996), no.~1, 149--179.


%
\bibitem[LX24]{LX24}
{\sc J.~Liu, D.~Xiang},
\emph{Exact global control of small divisors in rational normal form},
\href{https://doi.org/10.1088/1361-6544/ad4cd2}{Nonlinearity}, {\bf 37} 075020 






\bibitem[McC22]{McConnell}
{\sc R.~McConnell},
\emph{Nonlinear smoothing for the periodic generalized nonlinear Schr\"odinger equation,}
\href{https://doi.org/10.1016/j.jde.2022.09.017}{J. Differential Equations},   341 (2022), 353–379.

\bibitem[Pos96]{Pos96}
{\sc J.~P\"oschel},
\emph{A KAM theorem for some nonlinear PDEs},  
\href{http://www.numdam.org/item/ASNSP_1996_4_23_1_119_0}{Ann. Sc. Norm. Pisa Cl. Sci.}(4) \textbf{23} (1996), no.~1, 119--148.

%
%
\bibitem[Pos02]{Pos02}
{\sc J.~P{\"o}schel},
\emph{On the construction of almost-periodic solutions for a nonlinear Schr{\"o}dinger equation},
\href{https://doi.org/10.1017/S0143385702001086}{Ergodic Theory Dynam. Systems} \textbf{22} (2002), no.~5, 1537--1549.


\bibitem[PP15]{PP15}
{\sc C.~Procesi, M.~Procesi},
\emph{A KAM algorithm for the resonant non-linear Schr\"odinger equation}, 
\href{https://doi.org/10.1016/j.aim.2014.12.004}{Advances in Math.} (2015), 399--470.


%
%

\bibitem[Way90]{Way90} 
{\sc C.E.~Wayne}, 
\emph{Periodic and quasi-periodic solutions of nonlinear wave equations via KAM theory}, 
\href{https://doi.org/10.1007/BF02104499}{Comm. Math. Phys.} \textbf{127} (1990), no.~3, 479--528.


\end{thebibliography}
\end{document}